\documentclass[12pt]{amsart}
\usepackage{a4wide,tikz,hyperref}

\newtheorem{theorem}{Theorem}
\newtheorem{lemma}{Lemma}
\newtheorem{proposition}{Proposition}
\theoremstyle{definition}
\newtheorem{definition}{Definition}
\theoremstyle{remark}
\newtheorem{example}{Example}
\newtheorem{remark}{Remark}

\DeclareMathOperator\val{val}
\DeclareMathOperator\rep{rep}

\begin{document}

\begin{abstract}
Similarly to $\beta$-adic van der Corput sequences, abstract van der Corput sequences can be defined by abstract numeration systems.
Under some assumptions, these sequences are low discrepancy sequences.
The discrepancy function is computed explicitly, and the bounded remainder sets of the form $[0,y)$ are characterized. 
\end{abstract}

\title{Regularities of the distribution of abstract van der Corput sequences}
\author{Wolfgang Steiner}
\date\today
\address{LIAFA, CNRS, Universit\'e Paris Diderot -- Paris 7, case 7014, 75205 Paris Cedex 13, France}
\email{steiner@liafa.jussieu.fr}
\thanks{This work was supported by the French Agence Nationale de la Recherche, grant ANR--06--JCJC--0073.}

\maketitle

\section{Introduction}
Let $(x_n)_{n\ge0}$ be a sequence with $x_n\in[0,1)$ for all $n\ge 0$,
and 
\[
D(N,I) = \# \{0 \le n < N:\, x_n \in I\} - N \lambda(I)
\]
its \emph{discrepancy function} (or \emph{local discrepancy}) on the interval~$I$, where $\lambda(I)$ is the length of~$I$.
Then, $(x_n)_{n\ge0}$ is said to be a \emph{low discrepancy sequence} if $\sup_I D(N,I) = \mathcal{O}(\log N)$, where the supremum is taken over all intervals $I \subseteq [0,1)$. 
If $D(N,I)$ is bounded in~$N$, then $I$ is called a \emph{bounded remainder set}.
For details on discrepancy, we refer to \cite{Kuipers-Niederreiter74,Drmota-Tichy97}.
References to results on bounded remainder sets can be found in the introduction of \cite{Steiner06}.

In \cite{Barat-Grabner96,Ninomiya98a,Ninomiya98b}, $\beta$-adic van der Corput sequences are defined, and it is shown that they are low discrepancy sequences if $\beta$ is a Pisot number with irreducible $\beta$-polynomial.
Recall that a \emph{Pisot number} is an algebraic integer greater than~1 with all its conjugates lying in the interior of the unit disk. 
We refer to Section~\ref{beta} for the definition of the $\beta$-polynomial.
In \cite{Mori98,Ichikawa-mori04}, these results were extended to piecewise linear maps which generalize the $\beta$-transformation.
The proof in \cite{Ninomiya98a,Ninomiya98b} relies on the fact that cylinder sets of the $\beta$-transformation are bounded remainder sets if $\beta$ is a Pisot number with irreducible $\beta$-polynomial.
Under the same conditions on~$\beta$, bounded remainder sets of the form $[0,y)$, $0 \le y \le 1$, were completely characterized in \cite{Steiner06}: the $\beta$-expansion of~$y$ is finite or its tail is the same as that of the expansion of~1.

If $\beta$ is a Pisot number, then the language of the $\beta$-expansions is regular, which means that it is recognized by a finite automaton. 
Therefore, these $\beta$-expansions are special cases of abstract numeration systems as defined in \cite{Lecomte-Rigo01,Lecomte-Rigo02}, see Section~\ref{beta}.
In Section~\ref{ans}, we define van der Corput sequences related to more general abstract numeration systems.
Theorem~\ref{t1} in Section~\ref{disc} provides a new class of low discrepancy sequences. 
Finally, Theorem~\ref{t2} in Section~\ref{bounded} characterizes the bounded remainder sets of the form $[0,y)$ with respect to these abstract van der Corput sequences, generalizing the results in \cite{Steiner06}.

\section{Definitions and first results} \label{ans}

Let $(A,\le)$ be a finite and totally ordered alphabet. 
Denote by $A^*$ the free monoid generated by $A$ for the concatenation product, i.e., the set of finite words with letters in~$A$. 
The length of a word $w \in A^*$ is denoted by~$|w|$.
Extend the order on $A$ to $A^*$ by the \emph{shortlex} (or genealogical) \emph{order}, that is to say $v \le w$ if $v = w$ or $v < w$, where $v < w$ means that either $|v| < |w|$ or $|v| = |w|$ and there exist $p, v', w' \in A^*$, $a, b \in A$ such that $v = p a v'$, $w = p b w'$ and $a < b$. 

According to \cite{Lecomte-Rigo01}, the triple $S = (L,A,\le)$ is an \emph{abstract numeration system} if $L$ is an infinite regular language over~$A$, and the \emph{numerical value} of a word $w \in L$ is defined by  
\[
\val_S(w) = \# \{v \in L:\, v < w\}.
\]
If $\val_S(w)=n$, then we say that $w$ is the \emph{representation} of $n$ and write $\rep_S(n)=w$.

Denote by $A^\omega$ the set of (right) infinite words with letters in~$A$.
It is ordered by the \emph{lexicographical order}, that is to say $t \le u$ if $t = u$ or $t < u$, where $t < u$ means that there exist $p \in A^*$, $a, b \in A$, $t', u' \in A^\omega$, such that $t = p a t'$, $u = p b u'$ and $a < b$. 

Assume that the language $L$ grows exponentially, with
\[
\lim_{k\to\infty} \frac{\log\#\{v\in L:\,|v|\le k\}}{k} = \log\beta > 0.
\]
Suppose that $u \in A^\omega$ is the limit of words $w^{(k)} \in L$, i.e., every finite prefix of $u$ is a prefix of $w^{(k)}$ for all but a finite number of $k$'s.
Then, the value of $u$ is the real number
\begin{equation} \label{valomega}
\val_S^\omega(u) = \lim_{k\to\infty} \frac{\val_S(w^{(k)})}{\#\{v\in L:\,|v|\le|w^{(k)}|\}},
\end{equation}
if this limit exists and does not depend on the choice of $w^{(k)}$.
Conditions assuring the existence of this value are given in \cite{Lecomte-Rigo02}, see also Lemma~\ref{l3} and its proof. 
Let 
\[
L_\omega = \big\{u \in A^\omega:\, u = \lim_{k\to\infty} w^{(k)}\ \mbox{for some}\ w^{(k)} \in L\big\}.
\]
Since $\val_S^\omega(u) \in [1/\beta,1]$, we define the \emph{normalized value}
\[
\langle u\rangle = \frac{\beta\val_S^\omega(u)-1}{\beta-1} \in [0,1].
\]
We extend this definition to finite words $w \in L$ which are prefixes of words in $L_\omega$ by setting $\langle w\rangle = \langle u\rangle$, where $u$ is the smallest word in $L_\omega$ with prefix~$w$.
Since we want to define a sequence without multiple occurrences of the same value, we set
\[
L' = \{w \in L:\, \langle w\rangle \ne \langle v\rangle\ \mbox{for every}\ v\in L\ \mbox{with}\ v<w\}.
\]

Recall that the \emph{mirror image} of a word $w = w_1 w_2 \cdots w_k$, $w_j \in A$, is $\widetilde{w} = w_k \cdots w_2 w_1$ and that the mirror image of a language $L$ is $\widetilde{L} = \{\widetilde{w}:\, w \in L\}$.
Now, we are ready to define the main object of this paper, abstract van der Corput sequences.

\begin{definition}[Abstract van der Corput sequence]
Let $S = (L,A,\le)$ be an abstract numeration system, where $L$ is a regular language of exponential growth, every word $w \in L$ is the prefix of some infinite word $u \in L_\omega$, and the limit in (\ref{valomega}) exists for every $u \in L_\omega$.
Then, the \emph{abstract van der Corput sequence} corresponding to~$S$ is given by
\[
x_n = \langle w\rangle\quad \mbox{with}\quad \widetilde{w} = \rep_{\widetilde{S'}}(n),
\]
where $\widetilde{S'}$ is the abstract numeration system $\big(\widetilde{L'},A,\le\big)$.
\end{definition}

Thus, the set of values of an abstract van der Corput sequence is $\{x_n:\, n \ge 0\} = \{\langle w\rangle:\, w \in L\} = \{\langle w\rangle:\, w \in L'\}$, and the position of $\langle w\rangle$, $w \in L'$, in the sequence is determined by the shortlex order on the mirror image of~$L'$.
We need a number of further assumptions on the language~$L$ in order to get precise formulae for the discrepancy.
All these assumptions are satisfied by the $\beta$-adic van der Corput sequence when the language of the $\beta$-expansions is regular, cf.\ Section~\ref{beta}, and by Example~\ref{exa1} at the end of this section.

Let $\mathcal{A}_L = (Q,A,\tau,q_0,F)$ be a (complete) deterministic finite automaton recognizing~$L$, with set of states~$Q$, transition function $\tau:\, Q \times A \to Q$, initial state~$q_0$ and set of final states~$F$. 
The transition function is extended to words, $\tau:\, Q \times A^* \to Q$, by setting $\tau(q,\varepsilon) = q$ for the empty word $\varepsilon$ and $\tau(q,w a) = \tau(\tau(q,w),a)$.
A word $w \in A^*$ is accepted by $\mathcal{A}_L$, and thus in~$L$, if
and only if $\tau(q_0,w) \in F$. 

\begin{definition}[Totally ordered automaton]
A deterministic automaton $(Q,A,\tau,q_0,F)$ is said to be a \emph{totally ordered automaton} if there exists a total order on the set of states $Q$ such that, for all $q,r \in Q$,
\[
q \le r \quad \mbox{implies}\quad \tau(q,a) \le \tau(r,a)\ \mbox{for every}\ a \in A.
\]
\end{definition}

From now on, all automata will be totally ordered automata.
Furthermore, the maximal state will be the initial state and every state except the minimal one will be final. 
(In case $Q = F$, where the automaton recognizes $A^*$, we add a non-accessible state to~$Q$.)
W.l.o.g., the set of states will be $Q = \{0,1,\ldots,d\}$ for some positive integer~$d$, and the order on $Q$ will be the usual order on the integers, hence $q_0 = d$ and $F = \{1,\ldots,d\}$.
Moreover, we will assume that $\tau(0,a) = 0$ for every $a \in A$, i.e., the state $0$ is a sink.

\begin{lemma}\label{l1}
Let $L \subseteq A^*$ be recognized by a totally ordered automaton $\mathcal{A}_L = (Q, A, \tau, d, Q \setminus \{0\})$, with $Q = \{0,1,\ldots,d\}$ and $\tau(0,a) = 0$ for every $a \in A$.
Then, $\widetilde{L}$ is recognized by the totally ordered automaton $\mathcal{A}_{\widetilde{L}} = (Q, A, \widetilde{\tau}, 0, Q \setminus \{0\})$, where
\begin{equation} \label{etau}
\widetilde{\tau}(r,a) = \# \big\{q \in Q:\, \tau(q,a) + r > d\big\}\quad \mbox{for every}\ r \in Q,\, a \in A.
\end{equation}
In particular, we have $\widetilde{\tau}(0,a) = 0$ for every $a \in A$.
\end{lemma}

\begin{proof}
A deterministic automaton $\mathcal{A}' = (Q', A, \tau', q_0', F')$ recognizing $\widetilde{L}$ is obtained by determinizing the automaton which is given by inverting the transition function~$\tau$, see e.g.\ \cite{Sakarovitch09}.
This means that 
\begin{itemize}
\item
$Q'$ is a subset of the power set $\mathcal{P}(Q)$, 
\item
$\tau'(q',a) = \{q \in Q:\, \tau(q,a) \in q'\}$ for every $q' \in Q'$, $a \in A$,
\item
the set of final states in $\mathcal{A}_L$ is the initial state of $\mathcal{A}'$, i.e., $q_0' = \{1,\ldots,d\}$, 
\item
the final states in $\mathcal{A}'$ are those elements of $Q'$ which contain the initial state of~$\mathcal{A}_L$, i.e., $F' = \{q' \in Q':\, d \in q'\}$.
\end{itemize}
We show that $Q' \subseteq \{q_0', q_1', \ldots, q_d'\}$, where $q_r' = \{r+1,\ldots,d\}$ ($q_d'$ being the empty set). 
Since $\mathcal{A}_L$ is totally ordered and $\tau(0,a) = 0$, we obtain that 
\[
\tau'(q_0',a) = \{q \in Q:\, \tau(q,a) > 0\} = \{r+1,\ldots,d\} = q_r' \quad \mbox{for some}\ r \in Q.
\]
In the same way, we get, for every $r \in Q$ with $q_r' \in Q'$, that 
\[
\tau'(q_r',a) = \{q \in Q:\, \tau(q,a) > r\} = \{s+1,\ldots,d\} = q_s' \quad \mbox{for some}\ s \in Q.
\]
This shows that $Q' \subseteq \{q_0', q_1', \ldots, q_d'\}$.
It is easy to see that $\mathcal{A}'$ is a totally ordered automaton, with the order on $Q'$ given by $q_r' \le q_s'$ if $q_r' \subseteq q_s'$, i.e., $r \ge s$.
We clearly have $F' = \{q_0', q_1', \ldots, q_{d-1}'\} \cap Q'$.
If we extend the set of states to $\{q_0', q_1', \ldots, q_d'\}$ (with possibly non-accessible states) and label the states by $d-r$ instead of $q_r'$, we obtain $\mathcal{A}_{\widetilde{L}}$.
Therefore, $\mathcal{A}_{\widetilde{L}}$ is a totally ordered automaton, with $Q$ ordered by the usual order on the integers.
\end{proof}

The next lemma provides a fundamental characterization of the words in
a language~$L$ recognized by a totally ordered automaton $\mathcal{A}_L = (Q, A, \tau, d, Q \setminus \{0\})$ with $Q = \{0,1,\ldots,d\}$ and $\tau(0,a) = 0$ for every $a \in A$.

\begin{lemma} \label{l2}
Let $L, \tau, \widetilde{\tau}$ be as in Lemma~\ref{l1}, $w_1 \cdots w_k \in A^*$, $0 \le j \le k$.
We have $w_1 \cdots w_k \in L$ if and only if $\tau(d,w_1 \cdots w_j) + \widetilde{\tau}(d,w_k \cdots w_{j+1}) > d$.
\end{lemma}

\begin{proof}
Let $0 \le j \le k$. 
With the notation of the proof of Lemma~\ref{l1}, $\widetilde{\tau}(d,w_k \cdots w_{j+1}) = d-r$ can be written as $\tau'(q_0',w_k \cdots w_{j+1}) = q_r' = \{r+1,\ldots,d\}$. 
In~$\mathcal{A}_L$, this means that $w_{j+1} \cdots w_k$ leads to a final state from the state $q$, i.e., $\tau(q,w_{j+1} \cdots w_k) > 0$, if and only if $q > r$.
Therefore, we have $\tau(d,w_1 \cdots w_k) = \tau(\tau(d,w_1 \cdots w_j), w_{j+1} \cdots w_k) > 0$ if and only if $\tau(d,w_1 \cdots w_j) > r$, i.e., $\tau(d,w_1 \cdots w_j) + \widetilde{\tau}(d,w_k \cdots w_{j+1}) > d$. 
\end{proof}

\begin{remark} \label{Ferrers}
If we consider $\tau(d,a) + \cdots + \tau(1,a)$ as a partition of an integer, then $\widetilde{\tau}(d,a) + \cdots + \widetilde{\tau}(1,a)$ is the conjugate partition, since $\widetilde{\tau}(d-r,a) = \# \big\{q \in Q:\, \tau(q,a) > r\big\}$.
E.g., if $(\tau(d,a), \ldots, \tau(1,a)) = (4,2,1,0)$, then $(\widetilde{\tau}(d,a), \ldots, \widetilde{\tau}(1,a)) = (3,2,1,1)$, and the corresponding Ferrers diagram is
\[
\begin{array}{r l}
& \begin{array}{p{2mm}p{2mm}p{2mm}p{2mm}}3&2&1&1\end{array} \\
4 & \begin{array}{|p{2mm}|p{2mm}|p{2mm}|p{2mm}|}\hline&&& \\
\hline\end{array} \\
2 & \begin{array}{|p{2mm}|p{2mm}|}&\\
\hline\end{array} \\
1 & \begin{array}{|p{2mm}|}\\
\hline\end{array} \\
0 & \begin{array}{|p{2mm}}$ $\end{array}
\end{array}.
\]
\end{remark}

Next, we characterize the values of the abstract van der Corput sequence, under the assumption that the incidence matrix of the co-accessible part of~$\mathcal{A}_L$ is primitive.
(A~state $q$ is co-accessible if $\tau(q,w) \in F$ for some $w \in A^*$.)

\begin{lemma} \label{l3}
Let $L$ be as in Lemma~\ref{l1} and assume that $M_L = (\# \{a \in A:\, \tau(q,a) = r\})_{1\le q,r\le d}$ is a primitive matrix.
Then, the normalized value exists for every $u \in L_\omega$ and is given by
\begin{equation} \label{nvalue}
\langle u\rangle = \sum_{j=1}^\infty \epsilon_j(u) \beta^{-j}\quad \mbox{with}\quad \epsilon_j(u) = \sum_{a<u_j} \eta_{\tau(d,u_1 \cdots u_{j-1}a)},
\end{equation}
where $\beta$ is the Perron-Frobenius eigenvalue of $M_L$, $(\eta_1, \ldots, \eta_d)^t$ is the corresponding right (column) eigenvector with $\eta_d = 1$, $\eta_0 = 0$, and $u = u_1 u_2 \cdots$ with $u_j \in A$ for all $j \ge 1$.
\end{lemma}

\begin{proof}
The definition and primitivity of the incidence matrix $M_L$ give
\begin{equation} \label{Lqk}
\# \{v \in A^k:\, \tau(q,v) > 0\} = (0,\ldots,1,\ldots,0) M_L^k (1,\ldots,1)^t = c \eta_q \beta^k + \mathcal{O}(\rho^k)
\end{equation}
with constants $c > 0$ and $\rho < \beta$ such that every eigenvalue $\alpha \neq \beta$ of $M_L$ satisfies $|\alpha| < \rho$.
Due to the assumptions on $\mathcal{A}_L$, we have $u_1 \cdots u_k \in L$ for every $k \ge 1$.
Similarly to \cite{Lecomte-Rigo01,Lecomte-Rigo02}, we split up
\[
\{v \in L:\, v < u_1 \cdots u_k\} = \{v \in L:\, |v| < k\} \cup \bigcup_{1\le j\le k} \bigcup_{a<u_j} \{u_1 \cdots u_{j-1} a w \in L:\, w \in A^{k-j}\}.
\]
Since $\# \{v \in L:\, |v| \le k\} = \sum_{j=0}^k (c \beta^j + \mathcal{O}(\rho^j)) = c \frac{\beta^{k+1}}{\beta-1} + \mathcal{O}(\max(1,\rho)^k)$, we obtain
\[
\frac{\val_S(u_1\cdots u_k)}{\#\{v\in L:\,|v|\le k\}} = \frac{1}{\beta} + \frac{\beta-1}{\beta} \sum_{j=1}^k \sum_{a<u_j} \eta_{\tau(d,u_1\cdots u_{j-1}a)} \beta^{-j} + \mathcal{O}\bigg(\frac{\max(1,\rho)^k}{\beta^k}\bigg).
\]
Therefore, we have 
\[
\val_S^\omega(u) = \lim_{k\to\infty} \frac{\val_S(u_1\cdots u_k)}{\#\{v\in L:\,|v|\le k\}} = \frac{1}{\beta} + \frac{\beta-1}{\beta} \sum_{j=1}^\infty \epsilon_j(u) \beta^{-j},
\] 
and $\langle u\rangle = (\beta \val_S^\omega(u) - 1)/(\beta-1)$ yields (\ref{nvalue}).
\end{proof}

As a last preparation for the study of the discrepancy of abstract van der Corput sequences, we consider the language~$L'$.

\begin{lemma} \label{l4}
Let $L$ be as in Lemma~\ref{l3} and assume that $\tau(q,a_0) > 0$ for every $q > 0$, where $a_0$ denotes the smallest letter of~$A$.
Then, $L'$~consists exactly of those words in $L$ which do not end with~$a_0$.
\end{lemma}

\begin{proof}
We clearly have $\langle t\rangle \le \langle u\rangle$ if $t < u$, $t, u \in L_\omega$.
The primitivity of $M_L$ implies that $\eta_q > 0$ for all $q > 0$. Therefore, $\langle t\rangle = \langle u\rangle$ if and only if there exists no word $u' \in L_\omega$ with $t < u' < u$.
For $v < w$, $v, w \in L$, this means that $\langle v\rangle = \langle w\rangle$ if and only if $v$~is a prefix of~$w$ and $w$~is the smallest right extension of~$v$ in~$L$ of length~$|w|$.
Since $\tau(q,a_0) > 0$ for every $q > 0$, we have $v a_0^k \in L$ for all $k \ge 0$.
It follows that $\langle v\rangle = \langle w\rangle$ with $v < w$, $v, w \in L$, if and only if  $w = v a_0^{|w|-|v|}$.
Thus, $w \not\in L'$ if and only if $w$ ends with $a_0$.
\end{proof}

\begin{figure}[ht]
\includegraphics{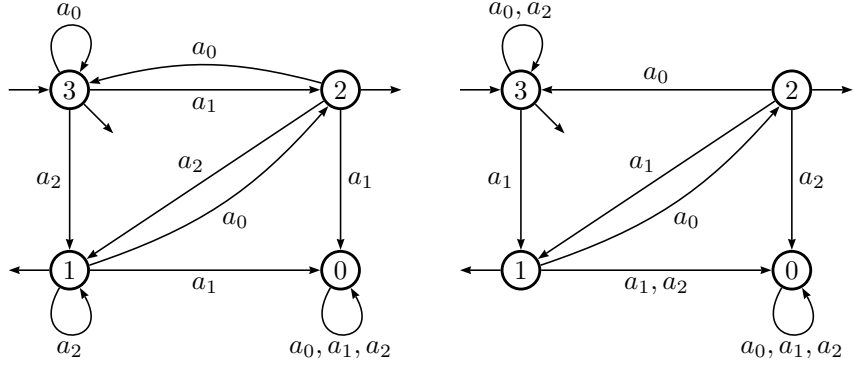}
\caption{The totally ordered automata $\mathcal{A}_L$ (left) and $\mathcal{A}_{\widetilde{L}}$ (right) of Example~\ref{exa1}.} \label{figexa1}
\end{figure}

\begin{example} \label{exa1}
Let $\mathcal{A}_L = (\{0,1,2,3\}, \{a_0,a_1,a_2\}, \tau, 3, \{1,2,3\})$ be the totally ordered automaton in Figure~\ref{figexa1} on the left. 
The first words in~$L$ (in the shortlex order) are
\begin{gather*}
\varepsilon, a_0, a_1, a_2, a_0 a_0, a_0 a_1, a_0 a_2, a_1 a_0, a_1 a_2, a_2 a_0, a_2 a_2, 
a_0 a_0 a_0, a_0 a_0 a_1, a_0 a_0 a_2, a_0 a_1 a_0, a_0 a_1 a_2, \\ a_0 a_2 a_0, a_0 a_2 a_2, a_1 a_0 a_0, a_1 a_0 a_1, a_1 a_0 a_2, a_1 a_2 a_0, a_1 a_2 a_2, a_2 a_0 a_0, a_2 a_0 a_2, a_2 a_2 a_0, a_2 a_2 a_2.
\end{gather*}
The transition functions $\tau$, $\widetilde{\tau}$ and the incidence matrices $M_L$, $M_{\widetilde{L}}$ are given by 
\[
\begin{array}{r|c c c}\tau & a_0 & a_1 & a_2 \\ \hline 0 & 0 & 0 & 0 \\ 1 & 2 & 0 & 1 \\ 2 & 3 & 0 & 1 \\ 3 & 3 & 2 & 1\end{array}, \quad M_L = \begin{pmatrix}1 & 1 & 0 \\ 1 & 0 & 1 \\ 1 & 1 & 1\end{pmatrix}, \quad \begin{array}{r|c c c}\widetilde{\tau} & a_0 & a_1 & a_2 \\ \hline 
0 & 0 & 0 & 0 \\ 1 & 2 & 0 & 0 \\ 2 & 3 & 1 & 0 \\ 3 & 3 & 1 & 3\end{array},\quad M_{\widetilde{L}} = \begin{pmatrix}0 & 1 & 0 \\ 1 & 0 & 1 \\ 1 & 0 & 2\end{pmatrix}.
\]
Recall that $M_L = (\# \{a \in A:\, \tau(q,a) = r\})_{1\le q,r\le 3}$, $M_{\widetilde{L}} = (\# \{a \in A:\, \widetilde{\tau}(q,a) = r\})_{1\le q,r\le 3}$ and that $\widetilde{\tau}$ can be calculated using Remark~\ref{Ferrers}.
Thus, $\widetilde{L}$ is recognized by the totally ordered automaton $\mathcal{A}_{\widetilde{L}}$ in Figure~\ref{figexa1} on the right.
The characteristic polynomial of $M_L$ is $x^3 - 2 x^2 - x + 1$, the dominant eigenvalue is $\beta \approx 2.247$, and $(\eta_1, \eta_2, \eta_3)^t = (\beta^2 - 2 \beta, -\beta^2 + 3 \beta - 1,1)^t \approx (0.555, 0.692, 1)^t$ is a right eigenvector of~$M_L$.
Since $x^3 - 2 x^2 - x + 1$ is irreducible, it must be the characteristic polynomial of $M_{\widetilde{L}}$ as well. 
The conditions of Lemma~\ref{l4} are satisfied, hence the first elements of the abstract van der Corput sequence corresponding to~$L$ are
\begin{gather*}
x_0 = \langle \varepsilon\rangle = 0,\
x_1 = \langle a_1\rangle = \frac{\eta_3}{\beta},\
x_2 = \langle a_2\rangle = \frac{\eta_3+\eta_2}{\beta},\
x_3 = \langle a_0 a_1\rangle = \frac{\eta_3}{\beta^2}, \\
x_4 = \langle a_0 a_2\rangle = \frac{\eta_3+\eta_2}{\beta^2},\,
x_5 = \langle a_1 a_2\rangle = \frac{\eta_3}{\beta} + \frac{\eta_3}{\beta^2},\
x_6 = \langle a_2 a_2\rangle = \frac{\eta_3+\eta_2}{\beta} + \frac{\eta_2}{\beta^2}, \\
x_7 = \langle a_0 a_0 a_1\rangle = \frac{\eta_3}{\beta^3},\
x_8 = \langle a_1 a_0 a_1\rangle = \frac{\eta_3}{\beta} + \frac{\eta_3}{\beta^3},\
x_9 = \langle a_0 a_0 a_2\rangle = \frac{\eta_3+\eta_2}{\beta^3}, \\
x_{10} = \langle a_1 a_0 a_2\rangle = \frac{\eta_3}{\beta} + \frac{\eta_3+\eta_2}{\beta^3},\
x_{11} = \langle a_2 a_0 a_2\rangle = \frac{\eta_3+\eta_2}{\beta} + \frac{\eta_3}{\beta^3}, \\
x_{12} = \langle a_0 a_1 a_2\rangle = \frac{\eta_3}{\beta^2} + \frac{\eta_3}{\beta^3},\
x_{13} = \langle a_0 a_2 a_2\rangle = \frac{\eta_3+\eta_2}{\beta^2} + \frac{\eta_2}{\beta^3}, \\
x_{14} = \langle a_1 a_2 a_2\rangle = \frac{\eta_3}{\beta} + \frac{\eta_3}{\beta^2} + \frac{\eta_2}{\beta^3},\
x_{15} = \langle a_2 a_2 a_2\rangle = \frac{\eta_3+\eta_2}{\beta} + \frac{\eta_2}{\beta^2} + \frac{\eta_2}{\beta^3}.
\end{gather*}
\centerline{\includegraphics{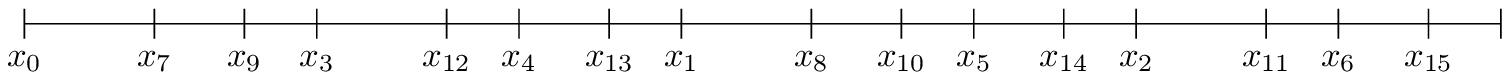}}
\end{example}

\section{$\beta$-adic van der Corput sequences} \label{beta}
To obtain Ninomiya's $\beta$-adic van der Corput sequences, consider a totally ordered automaton $\mathcal{A}_L = (\{0,1,\ldots,d\}, A, \tau, d, \{1,\ldots,d\})$ on the alphabet $A = \{0,1,\ldots,B\}$, with integers $b_q \in A$, $1 \le q \le d$, such that $\tau(q,a) = d$ for all $a < b_q$, and $\tau(q,a) = 0$ if and only if $a > b_q$ or $q = 0$.
Assume that $M_L$ is primitive (which is the only interesting case since $\tau(q,0) \ne d$ implies $b_q = 0$), and let $\beta$ be its Perron-Frobenius eigenvalue.
Then, we have
\[
\epsilon_j(u) = \sum_{a<u_j} \eta_{\tau(q,u_1 \cdots u_{j-1}a)} = \sum_{a<u_j} \eta_d = u_j \quad \mbox{for every}\ u = u_1 u_2 \cdots \in L_\omega,\ j \ge 1.
\]

Let $t_1 t_2 \cdots$ be the maximal sequence in $L_\omega$, i.e., $t_j = b_{\tau(d,t_1 \cdots t_{j-1})}$ for all $j \ge 1$.
Since $\mathcal{A}_L$ is a totally ordered automaton, $\tau(d,t_1 \cdots t_{j-1}) \le \tau(d,t_1 \cdots t_{k-1})$ implies $t_j < t_k$ or $t_j = t_k$, $\tau(d,t_1 \cdots t_j) \le \tau(d,t_1 \cdots t_k)$, thus $t_j t_{j+1} \cdots \le t_k t_{k+1} \cdots$, in particular $t_j t_{j+1} \cdots \le t_1 t_2 \cdots$.
By the maximality of $t_1 t_2 \cdots$, we have $\sum_{j=1}^\infty t_j \beta^{-j} = 1$.
The sequence $t_1 t_2 \cdots$ is called the \emph{infinite expansion of~$1$} in base~$\beta$ or \emph{quasi-greedy expansion} of~$1$, cf.\ \cite{Parry60}.
By the special structure of $\mathcal{A}_L$ and the primitivity of~$M_L$, we have $\{\tau(d,t_1 \cdots t_j):\, 0 \le j < d\} = \{1,\ldots,d\}$, hence $\tau(d,t_1 \cdots t_d) = \tau(d,t_1 \cdots t_m)$ for some $m < d$. 
Therefore, $t_1 t_2 \cdots$ is eventually periodic with preperiod length $m$ and period length $d-m$, which implies
\[
1 - \sum_{j=1}^d t_j \beta^{-j} = \bigg(1 - \sum_{j=1}^m t_j \beta^{-j}\bigg) \beta^{m-d},
\]
thus $\beta$ is a root of the polynomial 
\begin{equation} \label{betapol}
(x^d - t_1 x^{d-1} - \cdots - t_d x^0) - (x^m - t_1 x^{m-1} - \cdots - t_m x^0).
\end{equation}
If $\mathcal{A}_L$ is the minimal deterministic automaton recognizing~$L$ or, equivalently, $m$~and $d-m$ are the minimal preperiod and period lengths of $t_1 t_2 \cdots$, then (\ref{betapol}) is called \emph{$\beta$-polynomial}.

If $u = u_1 u_2 \cdots \in L_\omega$, then we have either $u = t_1 t_2 \cdots$ or $u_1 \cdots u_{k-1} = t_1 \cdots t_{k-1}$, $u_k < t_k$ for some $k \ge 1$.
Since $\tau(d,u_1 \cdots u_k) = d$ in the latter case, $u \in L_\omega$ is equivalent with $u_j u_{j+1} \cdots \le t_1 t_2 \cdots$ for all $j \ge 1$, cf.\ \cite{Parry60}. 
Therefore, $u$ is either the greedy or the quasi-greedy $\beta$-expansion of~$\langle u\rangle$.
Since $w 0 0 \cdots$ is a greedy $\beta$-expansion for every $w \in L$, the abstract van der Corput sequence given by $S = (L, A, \le)$ is exactly the $\beta$-adic van der Corput sequence defined in \cite{Ninomiya98a}.
Consequently, we call $\mathcal{A}_L$ a \emph{$\beta$-automaton}.

Conversely, let $t_1 t_2 \cdots$ be the infinite expansion of~$1$ in base $\beta > 1$, i.e., the unique sequence of integers satisfying $\sum_{j=1}^\infty t_j \beta^{-j} = 1$ and $0 0 \cdots < t_j t_{j+1} \cdots \le t_1 t_2 \cdots$ for all $j \ge 1$, cf.\ \cite{Parry60}.
Assume that $t_1 t_2 \cdots$ is eventually periodic.
In particular, this holds when $\beta$ is a Pisot number, see \cite{Bertrand77,Schmidt80}.
Let $d$ be the sum of the minimal preperiod and period lengths, and  $q_j = \# \{1 \le k \le d:\, t_k t_{k+1} \cdots \le t_{j+1} t_{j+2} \cdots\}$ for $j \ge 0$.
Then, the $\beta$-automaton is given by $b_{q_j} = t_{j+1}$ and $\tau(q_j,t_{j+1}) = q_{j+1}$, see Example~\ref{exa2}.

Note that $\beta$-adic van der Corput sequences were first considered in \cite{Barat-Grabner96} for some cases in which $\widetilde{L} = L$. 
In this case, the definition of the sequence is simpler beacause the $\beta$-expansion of $x_n$ is the mirror image of the expansion of~$n$ in a numeration system with respect to the linear recurrence corresponding to the $\beta$-polynomial.
In our notation, $\widetilde{L} = L$ is equivalent with $\widetilde{\tau} = \tau$, see also \cite{Brown-Yin00,Kwon09} for a different characterization.

\begin{example} \label{exa2}
Let $\beta$ be the real root of $x^3 - 4 x^2 - 2$.
Then, $t_1 t_2 \cdots = 4 0 1 4 0 1 \cdots$, hence we have $d = 3$, $q_0 = 3$, $q_1 = 1$, $q_2 = 2$, $q_3 = 3$, thus $b_3 = 4$, $b_1 = 0$, $b_2 = 1$, and  
\[
\begin{array}{r|c c c c c}\tau & 0 & 1 & 2 & 3 & 4 \\ \hline 0 & 0 & 0 & 0 & 0 & 0 \\ 1 & 2 & 0 & 0 & 0 & 0 \\ 2 & 3 & 3 & 0 & 0 & 0 \\ 3 & 3 & 3 & 3 & 3 & 1\end{array},\quad M_L = \begin{pmatrix}0 & 1 & 0 \\ 0 & 0 & 2 \\ 1 & 0 & 4\end{pmatrix}, \quad \begin{array}{r|c c c c c}\widetilde{\tau} & 0 & 1 & 2 & 3 & 4 \\ \hline 0 & 0 & 0 & 0 & 0 & 0 \\ 1 & 2 & 2 & 1 & 1 & 0 \\ 2 & 3 & 2 & 1 & 1 & 0 \\ 3 & 3 & 2 & 1 & 1 & 1\end{array}, \quad M_{\widetilde{L}} = \begin{pmatrix}2 & 2 & 0 \\ 2 & 1 & 1 \\ 3 & 1 & 1\end{pmatrix}.
\]
This example is also given in Section~2.3 in \cite{Steiner06}, with different notation. 
The substitution $\tau$ in \cite{Steiner06} plays the role of the transition function $\widetilde{\tau}$ in this paper.
\end{example}

\section{Discrepancy function} \label{disc}
For a given abstract van der Corput sequence $(x_n)_{n\ge0}$ and $y \in [0,1]$, we study now the behavior of the function $D(N,[0,y)) = \# \{0 \le n < N:\, x_n < y\} - N y$ as $N \to \infty$. 
Since $D(N,[y,z)) = D(N,[0,z)) - D(N,[0,y))$, this determines the discrepancy of $(x_n)_{n\ge0}$.

If $L$ satisfies the conditions of Lemma~\ref{l3}, then every $y \in [0,1]$ is the numerical value of some $u \in L_\omega$, see  \cite{Lecomte-Rigo02}.
We have the following lemma, where $a_0^\omega = a_0 a_0 \cdots$. 

\begin{lemma} \label{l5}
Let $(x_n)_{n\ge0}$ be an abstract van der Corput sequence with $L$ as in Lemma~\ref{l4}.
For $y = \langle u\rangle$ with $u = u_1 u_2 \cdots \in L_\omega$, $N \ge 0$ with $\rep_{\widetilde{S'}}(N) = w_\ell \cdots w_1$, we have
\[
\# \{0 \le n < N:\, x_n < y\} = \sum_{j=1}^\ell \sum_{a<u_j} \sum_{k=j+1}^\ell \sum_{b<w_k} \# L_{\tau(d,u_1\cdots u_{j-1}a),\widetilde{\tau}(d,w_\ell \cdots w_{k+1}b)}^{k-j-1} + C(N,u),
\]
where
$L_{q,r}^m = \{v \in A^m:\, \tau(q,v) + r > d\} = \{v \in A^m:\, q + \widetilde{\tau}(r,\widetilde{v}) > d\}$ and
\[
C(N,u) = \sum_{k=1}^\ell \# \{b < w_k:\, b w_{k+1} \cdots w_\ell a_0^\omega < u_k u_{k+1} \cdots,\, u_1 \cdots u_{k-1} b w_{k+1} \cdots w_\ell \in L\}.
\]
Since $\ell = \mathcal{O}(\log N)$, we have $C(N,u) = \mathcal{O}(\log N)$. 
\end{lemma}

\begin{proof}
We have to count the number of words $v \in L'$ with $\val_{\widetilde{S'}}(\widetilde{v}) < N$, i.e., $\widetilde{v} < w_\ell \cdots w_1$, and $\langle v\rangle < y$.
As in the proof of Lemma~\ref{l4}, we have $v a_0^{\ell-|v|} \in L$, thus $\langle v\rangle = \langle v a_0^{\ell-|v|}\rangle$. 
Since $w_\ell > a_0$, $\widetilde{v} < w_\ell \cdots w_1$ holds if and only if $a_0^{\ell-|v|} \widetilde{v} < w_\ell \cdots w_1$.
Therefore, we can count the number of words $v \in L \cap A^\ell$ with $\widetilde{v} < w_\ell \cdots w_1$ and $\langle v\rangle < y$, instead of those in~$L'$.

The inequality $\langle v\rangle < y$ is equivalent with $v a_0^\omega < u$. 
Thus, we have to count the $v$ in~$L$ of the form $v = u_1 \cdots u_{j-1} a v_{j+1} \cdots v_{k-1} b w_{k+1} \cdots w_\ell$ with $a < u_j$, $b < w_k$, $1 \le j < k \le \ell$, or $v = u_1 \cdots u_{k-1} b w_{k+1} \cdots w_\ell$ with $b < w_k$, $b w_{k+1} \cdots w_\ell a_0^\omega < u_k u_{k+1} \cdots$, $1 \le k \le \ell$.

Lemma~\ref{l2} yields that $u_1 \cdots u_{j-1} a v_{j+1} \cdots v_{k-1} b w_{k+1} \cdots w_\ell \in L$ if and only if $v_{j+1} \cdots v_{k-1} \in L_{\tau(d,u_1\cdots u_{j-1}a),\widetilde{\tau}(d,w_\ell \cdots w_{k+1}b)}^{k-j-1}$, which provides the main part of the formula. 
The other words give $C(N,u)$.
We have $\ell = \mathcal{O}(\log N)$ since $L$ grows exponentially and the same holds for $L'$ by Lemma~\ref{l4}.
Since $A$ is finite, we obtain $C(N,u) = \mathcal{O}(\log N)$.
\end{proof}

Similarly to \cite{Ninomiya98a,Ninomiya98b,Steiner06}, we assume now that the characteristic polynomial of $M_L$ is irreducible.
Let $\beta_2, \ldots, \beta_d$ be the conjugates of the Perron-Frobenius eigenvalue $\beta_1=\beta$.
For any $z \in \mathbb{Q}(\beta)$, denote by $z^{(i)} \in \mathbb{Q}(\beta_i)$ the image of $z$ by the isomorphism mapping $\beta$ to~$\beta_i$.
Similarly to~(\ref{Lqk}), we have some constants $\theta_1, \ldots, \theta_d \in \mathbb{Q}(\beta)$ such that
\[
\# L_{q,r}^k = (0,\ldots,1,\ldots,0) M_L^k (0,\ldots,0,1,\ldots,1)^t = \sum_{i=1}^d \eta_q^{(i)} \theta_r^{(i)} \beta_i^k\quad \mbox{for}\ 1 \le q, r \le d.
\]
Note that $(\theta_1, \ldots, \theta_d)^t$ is a right eigenvector of 
$M_{\widetilde{L}}$.
Set $\theta_0 = 0$ and define
\begin{equation} \label{defgamma}
\gamma_k(N) = \sum_{b<w_k} \theta_{\widetilde{\tau}(d,w_\ell \cdots w_{k+1}b)} \quad \mbox{for}\ N = \val_{\widetilde{S'}}(w_\ell \cdots w_1),\ 1 \le k \le \ell,
\end{equation}
similarly to $\epsilon_j(u)$.
Then,
\begin{gather*}
\sum_{a<u_j} \sum_{b<w_k} \# L_{\tau(d,u_1\cdots u_{j-1}a),\widetilde{\tau}(d,w_\ell \cdots w_{k+1}b)}^{k-j-1} = \sum_{i=1}^d \epsilon_j^{(i)}(u) \gamma_k^{(i)}(N) \beta_i^{k-j-1}, \\
N = \sum_{k=1}^\ell \sum_{b<w_k} \# L_{d,\widetilde{\tau}(d,w_\ell\cdots w_{k+1}b)}^{k-1} = \sum_{k=1}^\ell \sum_{b<w_k}\sum_{i=1}^d \eta_d^{(i)} \theta_{\widetilde{\tau}(d,w_\ell\cdots w_{k+1}b)}^{(i)} \beta_i^{k-1} =  \sum_{i=1}^d \sum_{k=1}^\ell \gamma_k^{(i)}(N) \beta_i^{k-1}.
\end{gather*}
For $y = \langle u\rangle$, we have thus
\begin{align} 
D(N,[0,y)) & = \sum_{j=1}^\ell \sum_{a<u_j} \sum_{k=j+1}^\ell \sum_{b<w_k} \# L_{\tau(d,u_1\cdots u_{j-1}a),\widetilde{\tau}(d,w_\ell \cdots w_{k+1}b)}^{k-j-1} + C(N,u) - N \sum_{j=1}^\infty \epsilon_j(u) \beta^{-j} \nonumber \\
& = \sum_{j=1}^\infty \sum_{i=1}^d \Bigg(\sum_{k=j+1}^\ell \epsilon_j^{(i)}(u) \gamma_k^{(i)}(N) \beta_i^{k-j-1} - \sum_{k=1}^\ell \gamma_k^{(i)}(N) \beta_i^{k-1} \epsilon_j(u) \beta^{-j}\Bigg) + C(N,u)\nonumber\\
& = \sum_{j=1}^\infty \Bigg(\sum_{k=j+1}^\ell \sum_{i=2}^d \gamma_k^{(i)}(N) \beta_i^{k-1} \Big(\epsilon_j^{(i)}(u) \beta_i^{-j} - \epsilon_j(u) \beta^{-j}\Big) \label{esumj} \\
& \hspace{15mm} - \sum_{k=1}^{\min(j,\ell)} \sum_{i=1}^d \gamma_k^{(i)}(N) \beta_i^{k-1} \epsilon_j(u) \beta^{-j}\Bigg) + C(N,u). \nonumber 
\end{align}
Since the series converges absolutely, we can change the order of summation, and get
\begin{align}
D(N,[0,y)) & = \sum_{k=1}^\ell \Bigg(\sum_{j=1}^{k-1} \sum_{i=2}^d \gamma_k^{(i)}(N) \beta_i^{k-1} \Big(\epsilon_j^{(i)}(u) \beta_i^{-j} - \epsilon_j(u) \beta^{-j}\Big) \nonumber \\
& \hspace{15mm} - \sum_{j=k}^\infty \sum_{i=1}^d \gamma_k^{(i)}(N) \beta_i^{k-1} \epsilon_j(u) \beta^{-j}\Bigg) + C(N,u). \label{esumk}
\end{align}

The following theorem states that $D(N,[0,y)) = \mathcal{O}(\log N)$ if $|\beta_i| < 1$ for $2 \le i \le d$, i.e., if $\beta$ is a Pisot number. 
The conditions on $M_L$ are subsumed in the following definition.

\begin{definition}[Pisot automaton]
A~deterministic automaton is said to be a \emph{Pisot automaton} if the incidence matrix of its restriction to the states which are both accessible and co-accessible has one simple eigenvalue $\beta > 1$, and all other eigenvalues satisfy $|\alpha| < 1$.
\end{definition}

\begin{theorem} \label{t1}
Let $S = (L,A,\le)$ be an abstract numeration system where $L$ is
recognized by a totally ordered Pisot automaton $\mathcal{A}_L = (\{0,1,\ldots,d\}, A, \tau, d, \{1,\ldots,d\})$ with $\tau(0,a) = 0$ for every $a \in A$, $\tau(q,a_0) > 0$ for the minimal letter $a_0 \in A$ and every $q > 0$.
Then, the corresponding abstract van der Corput sequence is a low
discrepancy sequence. 
\end{theorem}

\begin{proof}
If $\mathcal{A}_L$ is a Pisot automaton, then the characteristic polynomial of $M_L$ is irreducible, thus $M_L$ is primitive. 
Therefore, the conditions of Lemma~\ref{l5} are satisfied, and $D(N,[0,y))$ is given by~(\ref{esumk}).
Since $\gamma_k(N)$ and $\epsilon_j(u)$ take only finitely many different values and $\ell = \mathcal{O}(\log N)$, we obtain
\[
D(N,[0,y)) = \sum_{k=1}^\ell \mathcal{O}(1) + \mathcal{O}(\log N) = \mathcal{O}(\log N),
\]
where the constants implied by the $\mathcal{O}$-symbols do not depend on~$y$.
With $D(N,[y,z)) = D(N,[0,z)) - D(N,[0,y))$, we get $\sup_I D(N,I) = \mathcal{O}(\log N)$, and the theorem is proved.
\end{proof}

As Section~\ref{beta} shows, the $\beta$-adic van der Corput sequences defined by Pisot numbers~$\beta$ with irreducible $\beta$-polynomial considered in \cite{Ninomiya98a,Ninomiya98b} are special cases of the abstract van der Corput sequences in Theorem~\ref{t1}. 
By Lemma~\ref{l2}, $\widetilde{S} = (\widetilde{L}, A, \le)$ is another, usually different, abstract numeration system satisfying the conditions of Theorem~\ref{t1} whenever $\widetilde{\tau}(q,a_0) > 0$ for every $q > 0$, in particular when $\mathcal{A}_L$ is a $\beta$-automaton. 
The abstract van der Corput sequence considered in Example~\ref{exa1} is also a new low discrepancy sequence.

\section{Bounded remainder sets} \label{bounded}

Under the conditions of Theorem~\ref{t1}, (\ref{esumj}) gives
\begin{multline} \label{esumj2}
D(N,[0,y)) \\ = \sum_{j=1}^\infty \Bigg(\sum_{k=j+1}^\ell \sum_{i=2}^d \gamma_k^{(i)}(N) \epsilon_j^{(i)}(u) \beta_i^{k-j-1} - \sum_{k=1}^{\min(j,\ell)} \gamma_k(N) \epsilon_j(u) \beta^{k-j-1}\Bigg) + C(N,u) + \mathcal{O}(1)
\end{multline}
for $y = \langle u\rangle$.
If there exists some $m \ge 0$ such that $u_j = a_0$ for all $j > m$, then $\epsilon_j(u) = 0$ for all $j > m$ and $C(N,u)$ is bounded.
It follows that the interval $[0,\langle v\rangle) = [0,\langle v a_0^\omega\rangle)$ is a bounded remainder set for every finite word $v \in L$.

For any $\beta$-adic van der Corput sequence where $\beta$ is a Pisot number with irreducible $\beta$-polynomial, the bounded remainder sets $[0,y)$ have been characterized in \cite{Steiner06} by the fact that the tail of the $\beta$-expansion of $y$ is either $0^\omega$ or a suffix of $t_1 t_2 \cdots$. 
With our notation and the $\beta$-automaton defined in Section~\ref{beta}, this means that $\sum_{j=m+1}^\infty u_j \beta^{m-j} = \eta_q$ for some $m \ge 0$, $q \in Q$.
In the more general case, we have the following partial characterization. 
Note that the proof is simpler than the one in \cite{Steiner06}.

\begin{proposition} \label{p5}
Let $(x_n)_{n\ge0}$ be an abstract van der Corput sequence satisfying the conditions of Theorem~\ref{t1}, $u \in L_\omega$.
If there exists some $m \ge 0$, $q \in Q$, such that 
\[
\epsilon_{m+1}(u) \epsilon_{m+2}(u) \cdots = \epsilon_{q,1}(t_q) \epsilon_{q,2}(t_q) \cdots,
\] 
where $t_q = t_{q,1} t_{q,2} \cdots$ is the maximal sequence in $A^\omega$ with $\tau(q,t_{q,1} \cdots t_{q,j}) > 0$ for all $j \ge 1$ if $q > 0$, $t_0 = a_0^\omega$, $\epsilon_{q,j}(t_q) = \sum_{a<t_{q,j}} \eta_{\tau(q,t_{q,1} \cdots t_{q,j-1}a)}$, then $D(N,[0,\langle u\rangle))$ is bounded.
\end{proposition}

\begin{proof}
If $q = 0$, then $\epsilon_j(u) = 0$ for every $j > m$, and we have already seen that the boundedness of $D(N,[0,\langle u\rangle))$ follows from~(\ref{esumj2}).

For $q > 0$, note that $\epsilon_{q,2}(t_q) \epsilon_{q,3}(t_q) \cdots = \epsilon_{q',1}(t_{q'}) \epsilon_{q',2}(t_{q'}) \cdots$, where $q' = \tau(q,t_{q,1})$, and that $M_L$ is primitive.
Therefore, we can assume that $m$ is large enough such that $M_L^m$ has only positive entries. 
Hence, there exists some $v \in L \cap A^m$ such that $\tau(d,v) = q$.
Then, we have $v t_q \in L_\omega$ and
\[
\epsilon_{m+1}(v t_q) \epsilon_{m+2}(v t_q) \cdots = \epsilon_{q,1}(t_q) \epsilon_{q,2}(t_q) \cdots,
\]
thus $D\big(N,[0,\langle u\rangle)\big) = D\big(N,[0,\langle v t_q\rangle)\big) + \mathcal{O}(1)$ by~(\ref{esumj2}). 

It only remains to show that $D\big(N,[0,\langle v t_q\rangle)\big) = \mathcal{O}(1)$.
If $v t_q = t_d$, then $\langle v t_q\rangle = 1$ and $D\big(N,[0,\langle v t_q\rangle)\big) = 0$.
Otherwise, the successor $v'$ of $v$ in~$L$ has length~$m$, and $\langle v t_q\rangle = \langle v'\rangle$ since no sequence $u' \in L_\omega$ satisfies $v t_q < u' < v' a_0^\omega$, hence
\[
D\big(N,[0,\langle u\rangle)\big) = D\big(N,[0,\langle v t_q\rangle)\big) + \mathcal{O}(1) =  D\big(N,[0,\langle v'\rangle\big) + \mathcal{O}(1) = \mathcal{O}(1). \qedhere
\]
\end{proof}

Clearly $[y,z)$ is a bounded remainder set if both $[0,y)$ and $[0,z)$ have this property. 
For the converse, the following holds.

\begin{proposition} \label{p2}
Let $(x_n)_{n\ge0}$ be an abstract van der Corput sequence satisfying the conditions of Theorem~\ref{t1}. 
If $D(N,I)$ is bounded, then $\lambda(I) \in \mathbb{Q}(\beta)$.
\end{proposition}

\begin{proof}
The proof is very similar to that of Theorem~1 in \cite{Steiner06}. Therefore, we only give the main steps.
Define a substitution $q \mapsto \widetilde{\tau}(q,a_{q,0}) \cdots \widetilde{\tau}(q,a_{q,m_q})$, with $\{a_{q,0}, a_{q,1}, \ldots, a_{q,m_q}\} = \{a \in A:\, \widetilde{\tau}(q,a) > 0\}$ and $a_{q,0} < a_{q,1} < \cdots < a_{q,m_q}$, $1 \le q \le d$, which plays the role of the substitution $\tau$ in \cite{Steiner06}. 
Since $\widetilde{\tau}(d,a_0) = d$, we have $d \mapsto d w$ for some $w \in A^*$. 
Then, a continuous successor function on $L_\omega$ satisfying $\rep_{\widetilde{S'}}(n) a_0^\omega \mapsto \rep_{\widetilde{S'}}(n+1) a_0^\omega$ is topologically conjugate to the successor function $S$ on $\mathcal{D}$ defined in \cite{Steiner06}, see also \cite{Berthe-Rigo07}.
If $D(N,I)$ is bounded, then $\lambda(I)$ is an eigenvalue of the dynamical system $(\mathcal{D},S)$, see Theorem~5.1 in \cite{Shapiro78} and \cite{Steiner06}.
By Proposition~5 in \cite{Ferenczi-Mauduit-Nogueira96}, these eigenvalues are in $\mathbb{Q}(\beta)$.
\end{proof}

If $y \in \mathbb{Q}(\beta)$, then $y = \langle u\rangle$ for some eventually periodic sequence~$u$, see \cite{Rigo-Steiner05}. 
Let $p$ be the period length of $\epsilon_1(u) \epsilon_2(u) \cdots$ and $m$ the preperiod length.
From~(\ref{esumk}), we get that
\[
D(N,[0,y)) = \sum_{k=1}^\ell \!\Bigg(\!\sum_{i=2}^d \gamma_k^{(i)}(N) \sum_{j=1}^{k-1} \epsilon_j^{(i)}(u) \beta_i^{k-j-1} - \gamma_k(N) \sum_{j=k}^\infty \epsilon_j(u) \beta^{k-j-1}\!\Bigg)\! + C(N,u) + \mathcal{O}(1).
\]
Set $y_k = \sum_{j=k}^\infty \epsilon_j(u) \beta^{k-j-1}$.
For $k > m$, we have
\begin{align*}
y_k & = \frac{\epsilon_k(u)\beta^{p-1}+\cdots+\epsilon_{k+p-1}(u)}{\beta^p-1}, \\
\sum_{j=1}^{k-1} \epsilon_j^{(i)}(u) \beta_i^{k-j-1} & = \frac{\epsilon_{k-p}^{(i)}(u)\beta_i^{p-1}+\cdots+\epsilon_{k-1}^{(i)}(u)}{1-\beta_i^p} + \mathcal{O}(\beta_i^k) = -y_k^{(i)} + \mathcal{O}(\beta_i^k)\ \mbox{for}\ 2 \le i \le d,
\end{align*}
which gives
\[
D(N,[0,y)) = C(N,u) - \sum_{k=1}^\ell \sum_{i=1}^d \gamma_k^{(i)}(N) y_k^{(i)} + \mathcal{O}(1).
\]
Recall the definition of $\gamma_k(N)$ in (\ref{defgamma}), and set 
\begin{equation} \label{ezeta}
\zeta_r(z) = \sum_{i=1}^d \theta_r^{(i)} z^{(i)}\quad \mbox{for}\ r \in Q,\, z \in \mathbb{Q}(\beta).
\end{equation}
Then, we have
\begin{equation} \label{edelta}
D(N,[0,y)) = C(N,u) -  \sum_{k=1}^\ell \sum_{b<w_k} \zeta_{\widetilde{\tau}(d,w_\ell\cdots w_{k+1}b)}(y_k) + \mathcal{O}(1).
\end{equation}
Note that
\begin{equation} \label{ezetaeta}
\zeta_r(\eta_q) = \# L_{q,r}^0 = \left\{\begin{array}{cl}1 & \mbox{if}\ q + r > d, \\ 0 & \mbox{else.}\end{array}\right.
\end{equation}

The main result of this section is the following theorem.

\begin{theorem} \label{t2}
Let $(x_n)_{n\ge0}$ be an abstract van der Corput sequence defined by  an abstract numeration system $S = (L,A,\le)$, where $L$ is recognized by a totally ordered Pisot automaton $\mathcal{A}_L = (\{0,1,\ldots,d\}, A, \tau, d, \{1,\ldots,d\})$ with $\tau(0,a) = 0$ for every $a \in A$, $\tau(q,a_0) > q$ for the minimal letter $a_0 \in A$ and every $q \in \{1,\ldots,d-1\}$.
Then, $D(N,[0,\langle u\rangle))$, $u \in L_\omega$, is bounded in~$N$ if and only if $\langle u\rangle \in \mathbb{Q}(\beta)$ and there exists some $m \ge 0$ such that
\begin{equation} \label{ezeta2}
\zeta_{\widetilde{\tau}(d,\widetilde{v})}(y_k) = \left\{\begin{array}{cl}1 & \mbox{if}\ v a_0^\omega < u_k u_{k+1} \cdots\ \mbox{and}\ u_1 \cdots u_{k-1} v \in L, \\ 0 & \mbox{else,}\end{array}\right.
\end{equation}
for every $v \in A^*$, $k > m$, where $y_k = \sum_{j=k}^\infty \epsilon_j(u) \beta^{k-j-1}$, $\widetilde{\tau}$~is as in~(\ref{etau}), $\zeta$~as in~(\ref{ezeta}).
\end{theorem}

\begin{proof}
We show first that the conditions on $u \in L_\omega$ are sufficient for the boundedness of $D(N,[0,\langle u\rangle))$.
Note that the assumption $\tau(q,a_0) > q$ for $0 < q < d$ is not used here, but only the weaker assumption $\tau(q,a_0) > 0$ of Lemma~\ref{l4}.
For $N = \val_{\widetilde{S'}}(w_\ell \cdots w_1)$, we have
\begin{align*}
\sum_{k=1}^\ell \sum_{b<w_k} \zeta_{\widetilde{\tau}(d,w_\ell\cdots w_{k+1}b)}(y_k) & = \sum_{k=m+1}^\ell \# \bigg\{b < w_k: \begin{array}{l}b w_{k+1} \cdots w_\ell a_0^\omega < u_k u_{k+1} \cdots, \\ u_1 \cdots u_{k-1} b w_{k+1} \cdots w_\ell \in L\end{array}\bigg\} + \mathcal{O}(1) \\
& = C(N,u) + \mathcal{O}(1),
\end{align*}
thus $D(N,[0,\langle u\rangle)) = \mathcal{O}(1)$ by~(\ref{edelta}). 

\smallskip
Now, suppose that $D(N,[0,\langle u\rangle)) = \mathcal{O}(1)$ for $u \in L_\omega$.
By Proposition~\ref{p2}, we have $\langle u\rangle \in \mathbb{Q}(\beta)$, thus $u$ is eventually periodic by \cite{Rigo-Steiner05}.
Let $m \ge 0$, $p \ge 1$ be such that $u = u_1 \cdots u_m (u_{m+1} \cdots u_{m+p})^\omega$ and $\tau(d,u_1 \cdots u_{m}) = \tau(d,u_1 \cdots u_{m+p})$.

Consider $v \in A^*$ and $k > m$.
If $v \not\in L$, then (\ref{ezeta2}) holds since $\zeta_{\widetilde{\tau}(d,v)}(y_k) = \zeta_0(y_k) = 0$.
If $v = v' a_0$, then (\ref{ezeta2}) holds for $v$ if and only if it holds for~$v'$.
Therefore, we can assume $v \in L'$.
If $v$ is not the empty word, then we define integers
\[
N_h = \val_{\widetilde{S'}}\big((\widetilde{v} a_0^g)^h \widetilde{v} a_0^{k-1}\big) \quad \mbox{for}\ h \ge 0,
\]
with $g \ge \max(d-1,p)$ such that $g + |v|$ is a multiple of~$p$.
Set $v = v_k \cdots v_\ell$, with $k \le \ell$.

We show that
\begin{equation} \label{emu}
C(N_h,u) = (h+1) C(N_0,u)\quad \mbox{for all}\ h \ge 0.
\end{equation}
If $u_{m+1} \cdots u_{m+p} = a_0^p$, then $C(N_h,u) = 0$ and (\ref{emu}) holds.
If $u_{m+1} \cdots u_{m+p} > a_0^p$, then $g \ge p$ implies that $b v_{j+1} \cdots v_\ell a_0^\omega < u_j u_{j+1} \cdots$, $k \le j \le \ell$, if and only if $b v_{j+1} \cdots v_\ell (a_0^g v)^h a_0^\omega < u_j u_{j+1} \cdots$.
Since $\mathcal{A}_L$ is totally ordered, the assumption $\tau(q,a_0) > q$ for $0 < q < d$ implies $\tau(d,a_0) = d$ and $\tau(q,a_0^j) = d$ for all $q > 0$, $j \ge d-1$, in particular $\tau(q,a_0^g) = d$.
Therefore, we also have $u_1 \cdots u_{j-1} b v_{j+1} \cdots v_\ell \in L$ if and only if $u_1 \cdots u_{j-1} b v_{j+1} \cdots v_\ell (a_0^g v)^h \in L$.
Since $g+|v|$ is a multiple of the period length of~$u$, (\ref{emu}) follows. 

Since $\tau(q,a_0^g) = d$ for all $q > 0$, we also have $\widetilde{\tau}(q,a_0^g) = d$ for all $q > 0$, and thus $\zeta_{\widetilde{\tau}(d,(\widetilde{v}a_0^g)^h v_\ell\cdots v_{j+1}b)}(y_j) = \zeta_{\widetilde{\tau}(d,v_\ell\cdots v_{j+1}b)}(y_j)$.
By (\ref{edelta}), we get
\[
D(N_h,[0,\langle u\rangle)) = (h+1) \bigg(C(N_0,u) - \sum_{j=k}^\ell \sum_{b<v_j}\zeta_{\widetilde{\tau}(d,v_\ell\cdots v_{j+1}b)}(y_j)\bigg) + \mathcal{O}(1).
\]
Therefore, $D(N,[0,\langle u\rangle)) = \mathcal{O}(1)$ implies that
\begin{equation} \label{eN}
\sum_{j=k}^\ell \sum_{b<v_j} \zeta_{\widetilde{\tau}(d,v_\ell \cdots v_{j+1}b)}(y_j) = C(N_0,u)\quad \mbox{for}\ N_0 = \val_{\widetilde{S'}}(v_\ell \cdots v_k a_0^{k-1}),\ k > m.
\end{equation}

We will show that (\ref{ezeta2}) holds, by considering (\ref{eN}) for the integer $N_0'$ defined by the successor of $v$ in~$\widetilde{L'}$. 
Let first $v$ be the empty word, and $N_0' = \val_{\widetilde{S'}}(v_k' a_0^{k-1})$ with $v_k'$ such that $\widetilde{\tau}(d,b) = 0$ for all $b \in A$ with $a_0 < b < v_k'$.
Then, (\ref{eN}) gives
\[
\zeta_d(y_k) = \zeta_{\widetilde{\tau}(d,a_0)}(y_k) = C(N_0',u) = \left\{\begin{array}{cl}1 & \mbox{if}\ a_0^\omega < u_k u_{k+1} \cdots, \\ 0 & \mbox{else,}\end{array}\right. 
\]
thus (\ref{ezeta2}) holds if $v$ is the empty word.
Assume next that there exists some $v_k' > v_k$ with $\widetilde{\tau}(d,v_\ell \cdots v_{k+1} v_k') > 0$.
Choose $v_k'$ such that $\widetilde{\tau}(d,v_\ell \cdots v_{k+1} b) = 0$ for all $b \in A$ with $v_k < b < v_k'$, and set $N_0' = \val_{\widetilde{S'}}(v_\ell \cdots v_{k+1} v_k' a_0^{k-1})$.
Then, (\ref{ezeta2}) holds because of
\begin{multline*}
\zeta_{\widetilde{\tau}(d,v_\ell\cdots v_k)}(y_k) = \sum_{b<v_k'} \zeta_{\widetilde{\tau}(d,v_\ell\cdots v_{k+1} b)}(y_k) - \sum_{b<v_k} \zeta_{\widetilde{\tau}(d,v_\ell\cdots v_{k+1}b)}(y_k) \\
= C(N_0',u) - C(N_0,u) = \left\{\begin{array}{cl}1 & \mbox{if}\ v_k \cdots v_\ell a_0^\omega < u_k u_{k+1} \cdots,\, u_1 \cdots u_{k-1} v_k \cdots v_\ell \in L, \\ 0 & \mbox{else.}\end{array}\right. 
\end{multline*}

Now, we proceed by induction. 
We know that (\ref{ezeta2}) holds if $v$ is the empty word.
Assume that (\ref{ezeta2}) holds for all words $v$ of length $\ell - k$, and consider
\begin{align}
\zeta_{\widetilde{\tau}(d,v_\ell \cdots v_{k+1})}(y_{k+1}) & = \sum_{i=1}^d \theta_{\widetilde{\tau}(d,v_\ell\cdots v_{k+1})}^{(i)}
\big(\beta y_k-\epsilon_k(u)\big)^{(i)} \nonumber \\
& = \sum_{i=1}^d \beta_i \theta_{\widetilde{\tau}(d,v_\ell\cdots v_{k+1})}^{(i)} y_k^{(i)} - \sum_{i=1}^d \theta_{\widetilde{\tau}(d,v_\ell\cdots v_{k+1})}^{(i)} \sum_{b<u_k} \eta_{\tau(d,u_1\cdots u_{k-1}b)}^{(i)}. \label{ezeta3}
\end{align}
Let $v_k$ be the maximal letter with $\widetilde{\tau}(d,v_\ell \cdots v_{k+1} v_k) > 0$. 
Since $(\theta_1, \ldots, \theta_d)^t$ is an eigenvector of~$M_{\widetilde{L}}$ and we already know that $\zeta_{\widetilde{\tau}(d,v_\ell \cdots v_{k+1} b)}(y_k)$, $b < v_k$, is given by (\ref{ezeta2}), we have 
\begin{multline*}
\sum_{i=1}^d \beta_i \theta_{\widetilde{\tau}(d,v_\ell\cdots v_{k+1})}^{(i)} y_k^{(i)} = \sum_{i=1}^d \sum_{b\in A} \theta_{\widetilde{\tau}(d,v_\ell\cdots v_{k+1}b)}^{(i)} y_k^{(i)} = \sum_{b\in A} \zeta_{\widetilde{\tau}(d,v_\ell\cdots v_{k+1}b)}(y_k) \\
= \zeta_{\widetilde{\tau}(d,v_\ell\cdots v_k)}(y_k) + \# \{b < v_k:\, b v_{k+1} \cdots v_\ell a_0^\omega < u_k u_{k+1} \cdots,\, u_1 \cdots u_{k-1} b v_{k+1} \cdots v_\ell \in L\}.
\end{multline*}
Using this equation,
\begin{align*}
\sum_{i=1}^d \theta_{\widetilde{\tau}(d,v_\ell\cdots v_{k+1})}^{(i)} \sum_{b<u_k} \eta_{\tau(d,u_1\cdots u_{k-1}b)}^{(i)} & = \sum_{b<u_k} \# L_{\tau(d,u_1\cdots u_{k-1}b),\widetilde{\tau}(d,v_\ell\cdots v_{k+1})}^0 \\
& = \# \{b < u_k:\, u_1 \cdots u_{k-1} b v_{k+1} \cdots v_\ell \in L\}
\end{align*}
and the induction hypothesis, (\ref{ezeta3}) yields
\begin{align*}
\zeta_{\widetilde{\tau}(d,v_\ell\cdots v_k)}(y_k) & = \zeta_{\widetilde{\tau}(d,v_\ell \cdots v_{k+1})}(y_{k+1}) + \# \{b < u_k:\, u_1 \cdots u_{k-1} b v_{k+1} \cdots v_\ell \in L\} \\
& \quad - \# \{b < v_k:\, b v_{k+1} \cdots v_\ell a_0^\omega < u_k u_{k+1} \cdots,\, u_1 \cdots u_{k-1} b v_{k+1} \cdots v_\ell \in L\} \\
& = \left\{\begin{array}{cl}1 & \mbox{if}\ v_{k+1} \cdots v_\ell a_0^\omega < u_{k+1} u_{k+2} \cdots\ \mbox{and}\ u_1 \cdots u_k v_{k+1} \cdots v_\ell \in L, \\ 0 & \mbox{else,}\end{array}\right. \\
& \quad + \left\{\begin{array}{cl}1 & \mbox{if}\ v_k < u_k\ \mbox{and}\ u_1 \cdots u_{k-1} v_k \cdots v_\ell \in L, \\ \! -1 & \mbox{if}\ v_k > u_k,\, v_{k+1} \cdots v_\ell a_0^\omega < u_{k+1} u_{k+2} \cdots,\, u_1 \cdots u_k v_{k+1} \cdots v_\ell \in L, \\ 0 & \mbox{else,}\end{array}\right. \\
& = \left\{\begin{array}{cl}1 & \mbox{if}\ v_k \cdots v_\ell a_0^\omega < u_k u_{k+1} \cdots\ \mbox{and}\ u_1 \cdots u_{k-1} v_k \cdots v_\ell \in L, \\ 0 & \mbox{else,}\end{array}\right.
\end{align*}
where we have used that $u_1 \cdots u_{k-1} b v_{k+1} \cdots v_\ell$ can be in $L$ only if $b \le v_k$.
Therefore, (\ref{ezeta2}) holds for all words $v$ of length $\ell - k + 1$, and Theorem~\ref{t2} is proved.
\end{proof}

In the case of $\beta$-adic van der Corput sequences, all bounded remainder sets satisfy the conditions of Proposition~\ref{p5}, see~\cite{Steiner06}.
The following example shows that this is probably not true in the more general setting, i.e., that there might be sequences $u$ satisfying the conditions of Theorem~\ref{t1}, but not those of Proposition~\ref{p5}.
However, $\mathcal{A}_L$ is not a Pisot automaton in this example, and we have not found an example where $\mathcal{A}_L$ satisfies the conditions of Theorem~\ref{t1}.

\begin{example}
Let $\mathcal{A}_L = (\{0,1,2,3,4\}, \{a_0,a_1,a_2,a_3\}, \tau, 4, \{1,2,3,4\})$ be the totally ordered automaton with $\tau, \widetilde{\tau}$ given by the transition tables
\[
\begin{array}{r|c c c c}\tau & a_0 & a_1 & a_2 & a_3 \\ \hline 0 & 0 & 0 & 0 & 0 \\ 1 & 3 & 2 & 1 & 0 \\ 2 & 3 & 2 & 2 & 1 \\ 3 & 4 & 4 & 3 & 1 \\ 4 & 4 & 4 & 4 & 2\end{array}, \quad \begin{array}{r|c c c c}\widetilde{\tau} & a_0 & a_1 & a_2 & a_3 \\ \hline 0 & 0 & 0 & 0 & 0 \\ 1 & 2 & 2 & 1 & 0 \\ 2 & 4 & 2 & 2 & 0 \\ 3 & 4 & 4 & 3 & 1 \\ 4 & 4 & 4 & 4 & 3\end{array}.
\]
If $u = a_3 a_0 a_2^\omega$, then $\tau(4,u_1 \cdots u_{k-1}) = 3$ and thus $\epsilon_k(u) = 2 \eta_4$ for all $k \ge 3$.
We obtain $y_k = \eta_3 - \eta_2 + \eta_1$, which implies $\zeta_4(y_k) = \zeta_2(y_k) = 1$, $\zeta_3(y_k) = \zeta_1(y_k) = 0$ by (\ref{ezetaeta}). 
It can be easily verified that (\ref{ezeta2}) holds for all $v \in A^*$, $k \ge 3$, but the conditions on~$u$ of Proposition~\ref{p5} are not satisfied.
\end{example}

We conclude by the remark that the boundedness of $D(N,I)$ is not invariant under translation of the interval, i.e., $D(N,[y,z))$ can be unbounded if $[0,z-y)$ is a bounded remainder set and vice versa, see \cite{Steiner06}.

\section*{Acknowledgements}
I am grateful to Philippe Nadeau for indicating the links with
partitions to me, and to the referee for many suggestions improving the quality of the paper. 

\bibliographystyle{amsalpha}
\bibliography{local2}

\providecommand{\bysame}{\leavevmode\hbox to3em{\hrulefill}\thinspace}
\providecommand{\MR}{\relax\ifhmode\unskip\space\fi MR }
\providecommand{\MRhref}[2]{%
  \href{http://www.ams.org/mathscinet-getitem?mr=#1}{#2}
}
\providecommand{\href}[2]{#2}
\begin{thebibliography}{FMN96}

\bibitem[Ber77]{Bertrand77}
A.~Bertrand, \emph{D\'eveloppements en base de {P}isot et r\'epartition modulo
  {$1$}}, C. R. Acad. Sci. Paris S\'er. A-B \textbf{285} (1977), no.~6,
  A419--A421.

\bibitem[BG96]{Barat-Grabner96}
G.~Barat and P.~J. Grabner, \emph{Distribution properties of {$G$}-additive
  functions}, J. Number Theory \textbf{60} (1996), no.~1, 103--123.

\bibitem[BR07]{Berthe-Rigo07}
V.~Berth{\'e} and M.~Rigo, \emph{Odometers on regular languages}, Theory
  Comput. Syst. \textbf{40} (2007), no.~1, 1--31.

\bibitem[BY00]{Brown-Yin00}
G.~Brown and Q.~Yin, \emph{{$\beta$}-transformation, natural extension and
  invariant measure}, Ergodic Theory Dynam. Systems \textbf{20} (2000), no.~5,
  1271--1285.

\bibitem[DT97]{Drmota-Tichy97}
M.~Drmota and R.~F. Tichy, \emph{Sequences, discrepancies and applications},
  Lecture Notes in Mathematics, vol. 1651, Springer-Verlag, Berlin, 1997.

\bibitem[FMN96]{Ferenczi-Mauduit-Nogueira96}
S.~Ferenczi, C.~Mauduit, and A.~Nogueira, \emph{Substitution dynamical systems:
  algebraic characterization of eigenvalues}, Ann. Sci. \'Ecole Norm. Sup. (4)
  \textbf{29} (1996), no.~4, 519--533.

\bibitem[IM04]{Ichikawa-mori04}
Y.~Ichikawa and M.~Mori, \emph{Discrepancy of van der {C}orput sequences
  generated by piecewise linear transformations}, Monte Carlo Methods Appl.
  \textbf{10} (2004), no.~2, 107--116.

\bibitem[KN74]{Kuipers-Niederreiter74}
L.~Kuipers and H.~Niederreiter, \emph{Uniform distribution of sequences},
  Wiley-Interscience [John Wiley \& Sons], New York, 1974, Pure and Applied
  Mathematics.

\bibitem[Kwo09]{Kwon09}
D.~Y. Kwon, \emph{The natural extensions of {$\beta$}-transformations which
  generalize baker's transformations}, Nonlinearity \textbf{22} (2009), no.~2,
  301--310.

\bibitem[LR01]{Lecomte-Rigo01}
P.~B.~A. Lecomte and M.~Rigo, \emph{Numeration systems on a regular language},
  Theory Comput. Syst. \textbf{34} (2001), no.~1, 27--44.

\bibitem[LR02]{Lecomte-Rigo02}
P.~Lecomte and M.~Rigo, \emph{On the representation of real numbers using
  regular languages}, Theory Comput. Syst. \textbf{35} (2002), no.~1, 13--38.

\bibitem[Mor98]{Mori98}
M.~Mori, \emph{Low discrepancy sequences generated by piecewise linear maps},
  Monte Carlo Methods Appl. \textbf{4} (1998), no.~2, 141--162.

\bibitem[Nin98a]{Ninomiya98a}
S.~Ninomiya, \emph{Constructing a new class of low-discrepancy sequences by
  using the {$\beta$}-adic transformation}, Math. Comput. Simulation
  \textbf{47} (1998), no.~2-5, 403--418.

\bibitem[Nin98b]{Ninomiya98b}
\bysame, \emph{On the discrepancy of the {$\beta$}-adic van der {C}orput
  sequence}, J. Math. Sci. Univ. Tokyo \textbf{5} (1998), no.~2, 345--366.

\bibitem[Par60]{Parry60}
W.~Parry, \emph{On the {$\beta $}-expansions of real numbers}, Acta Math. Acad.
  Sci. Hungar. \textbf{11} (1960), 401--416.

\bibitem[RS05]{Rigo-Steiner05}
M.~Rigo and W.~Steiner, \emph{Abstract {$\beta$}-expansions and ultimately
  periodic representations}, J. Th\'eor. Nombres Bordeaux \textbf{17} (2005),
  no.~1, 283--299.

\bibitem[Sak09]{Sakarovitch09}
J.~Sakarovitch, \emph{Elements of automata theory}, Cambridge University Press,
  Cambridge, 2009.

\bibitem[Sch80]{Schmidt80}
K.~Schmidt, \emph{On periodic expansions of {P}isot numbers and {S}alem
  numbers}, Bull. London Math. Soc. \textbf{12} (1980), no.~4, 269--278.

\bibitem[Sha78]{Shapiro78}
L.~Shapiro, \emph{Regularities of distribution}, Studies in probability and
  ergodic theory, Adv. in Math. Suppl. Stud., vol.~2, Academic Press, New York,
  1978, pp.~135--154.

\bibitem[Ste06]{Steiner06}
W.~Steiner, \emph{Regularities of the distribution of {$\beta$}-adic van der
  {C}orput sequences}, Monatsh. Math. \textbf{149} (2006), no.~1, 67--81.

\end{thebibliography}
\end{document}